\documentstyle[10pt]{article}

\font\tenof=msym10 
\def\su{\mbox{\rm su$_q$(2)}}
\def\C{\mbox{\tenof C}}
\def\R{\mbox{\tenof R}}
\def\N{\mbox{\tenof N}}
\def\alg{\mbox{${\cal A}^+_q(1)$}}
\def\case#1#2{{\textstyle{#1\over #2}}}
\def\cR{\mbox{${\cal R}$}}
\def\id{\mbox{\rm id}}

\begin{document}
%
%
\begin{center}
{\large\bf NONLINEAR DEFORMED su(2) ALGEBRAS INVOLVING TWO
DEFORMING FUNCTIONS \footnote{Presented at the 5th International Colloquium
``Quantum Groups and Integrable  Systems'', Prague, 20--22 June 1996}}

\bigskip\bigskip
{D. BONATSOS, P. KOLOKOTRONIS} 

\medskip
{\it Institute of Nuclear Physics, NCSR Demokritos, GR-15310 Aghia Paraskevi,
Attiki, Greece}

\bigskip
{C. DASKALOYANNIS}

\medskip
{\it Department of Physics, Aristotle University of Thessaloniki, GR-54006
Thessaloniki, Greece}

\bigskip
{A. LUDU}

\medskip
{\it Department of Theoretical Physics, Faculty of Physics, University of Bucharest,
Bucharest-Magurele, P.O.~Box MG-5211, Romania}

\bigskip
{C. QUESNE \footnote{Directeur de recherches FNRS}}

\medskip
{\it Physique Nucl\'eaire Th\'eorique et Physique Math\'ematique, Universit\'e
Libre de Bruxelles, Campus de la Plaine CP229, Boulevard du Triomphe, B-1050
Brussels, Belgium}

\end{center}
 
\begin{abstract}
The most common nonlinear deformations of the su(2) Lie algebra, introduced by
Polychronakos and Ro\v cek, involve a single arbitrary function of $J_0$ and
include the quantum algebra su$_q$(2) as a special case. In the present contribution,
less common nonlinear deformations of su(2), introduced by Delbecq and Quesne and
involving two deforming functions of $J_0$, are reviewed. Such algebras include
Witten's quadratic deformation of su(2) as a special case. Contrary to the former
deformations, for which the spectrum of $J_0$ is linear as for su(2), the latter give
rise to exponential spectra, a property that has aroused much interest in connection
with some physical problems. Another interesting algebra of this type, denoted by
\alg, has two series of ($N+1$)-dimensional unitary irreducible representations,
where $N=0$, 1, 2,~$\ldots$. To allow the coupling of any two such representations,
a generalization of the standard Hopf axioms is proposed. The resulting algebraic
structure, referred to as a two-colour quasitriangular Hopf algebra, is described.
\end{abstract} 
%
%
\section{Introduction}
Quantized universal enveloping algebras, also called $q$-algebras, refer to some
specific deformations of (the universal enveloping algebra of) Lie algebras, to
which they reduce when the deforming parameter~$q$ (or set of deforming
parameters) goes to one~\cite{drinfeld}. The simplest example of $q$-algebra,
\su\ $\equiv$ U$_q$(su(2)), was first introduced by Sklyanin, and by Kulish and
Reshetikhin~\cite{sklyanin}. It has found a lot of applications in various branches of
physics since its realization in terms of $q$-bosonic operators was proposed by
Biedenharn and Macfarlane~\cite{biedenharn}.\par
%
%
The \su\ algebra is a special case of more general deformations of su(2),
independently introduced by Polychronakos and Ro\v cek~\cite{poly}.
They involve one arbitrary function $f(J_0)$ in the commutator of~$J_+$
with~$J_-$, and their representation theory is characterized by a rich variety of
phenomena, whose interest in particle physics has been stressed.\par
%
%
More recently, deformations of su(2) involving two deforming functions~$F(J_0)$
and~$G(J_0)$ in the commutator of~$J_+$ with~$J_-$ and in that of~$J_0$
with~$J_+$ or~$J_-$, respectively, have been proposed by Delbecq and
Quesne~\cite{delbecq1,delbecq2,delbecq3}. It is the purpose of the present
contribution to review the construction and representation theory of such algebras,
and to show how the problem of endowing some of them with a Hopf algebraic
structure can be addressed~\cite{bonatsos}.\par
%
%
\section{Nonlinear deformed su(2) algebras}
Polychronakos-Ro\v cek algebras (PRA's) are associative algebras over $\C$,
generated by three operators
$j_0 = (j_0)^{\dagger}$, $j_+$, and $j_- = (j_+)^{\dagger}$, satisfying the
commutation relations \cite{poly}
\begin{equation}
[j_0, j_+] = j_+, \qquad [j_0, j_-] = - j_-, \qquad [j_+, j_-] = f(j_0),
\label{eq:PRA}
\end{equation} 
where $f(z)$ is a real, parameter-dependent function of $z$, holomorphic in
the neighbourhood of zero, and going to $2z$ for some values of the
parameters. These algebras have a Casimir operator given by $c= j_- j_+ + h(j_0) =
j_+ j_- + h(j_0) - f(j_0)$, in terms of another real function $h(z)$, related to $f(z)$
through the equation $h(z) - h(z - 1) = f(z)$. An explicit expression for $h(z)$ has
been given by Delbecq and Quesne~\cite{delbecq1} in terms of Bernoulli polynomials
and Bernoulli numbers.\par
%
%
{}For all PRA's, the spectrum of~$j_0$ is linear as in the special case of the
$q$-algebra \su. The latter corresponds to $f(j_0) = [2j_0]_q$, and $h(j_0) = [j_0]_q
[j_0 + 1]_q$, with $[x]_q \equiv (q^x - q^{-x})/(q - q^{-1})$,  and  $q\in\R^+$ (the case
where $q$ is a phase will not be considered here, as throughout the present work we
shall restrict the parameters to real values)~\cite{biedenharn}.\par
%
%
Delbecq-Quesne algebras (DQA's) differ from PRA's by the replacement
of~(\ref{eq:PRA}) by~\cite{delbecq1}
\begin{equation}
  [J_0, J_+] = G(J_0) J_+, \qquad [J_0, J_-] = - J_- G(J_0), \qquad
  [J_+, J_-] = F(J_0),       \label{eq:DQA}
\end{equation}
where $J_0 = (J_0)^{\dagger}$, $J_- = (J_+)^{\dagger}$, and the commutators involve
two real, parameter-dependent functions of $z$, $F(z)$ and $G(z)$, holomorphic in
the neighbourhood of zero, and going to $2z$ and~1 for some values of the
parameters, respectively. These functions are further restricted by the assumption
that the algebras have a Casimir operator given by $C = J_- J_+ + H(J_0) = J_+ J_- +
H(J_0) - F(J_0)$, in terms of some real function $H(z)$, holomorphic in the
neigh\-bour\-hood  of zero. The latter restriction implies that $F(z)$, $G(z)$, and
$H(z)$ satisfy the consistency condition $H(z) - H\bigl(z - G(z)\bigr) = F(z)$.\par
%
%
Since for $G(J_0)=1$, DQA's reduce to PRA's, the first significant case corresponds
to $G(J_0) = 1 + (1-q) J_0$, where $q \in \R^+$. In such a case, it has been shown that
there exist ($\lambda-1$)-parameter algebras ${\cal A}^+_{\alpha_2 \alpha_3
\ldots \alpha_{\lambda-1} q}(\lambda,1)$, for which the functions $F(J_0)$ and
$H(J_0)$ are polynomials of degree~$\lambda$ in $J_0$~\cite{delbecq1}. In
particular, for $\lambda=2$ and~3, one finds the algebras ${\cal
A}^+_q(2,1)$~\cite{delbecq2} and ${\cal A}^+_{p,q}(3,1)$~\cite{delbecq3}, for which 
$F(J_0) = 2J_0 \bigl(1 + (1-q) J_0\bigr)$, $H(J_0) = 2 (1+q)^{-1} J_0 (J_0 + 1)$, and
$F(J_0) = 2J_0 \bigl(1 + (1-q) J_0\bigr) \bigl(1 - (1-p) J_0\bigr)$, $H(J_0) = 2
\bigl((1+q)(1+q+q^2)\bigr)^{-1} J_0 (J_0+1) \bigl(1 + (p+q)q - (1-p)(1+q)J_0\bigr)$,
respectively. The former is equivalent to Witten's first deformation of
su(2)~\cite{witten}.\par
%
%
The representation theory of the DQA's can be dealt with as that of \su, or more
generally of the PRA's. Considering the case where $G(J_0)= 1 + (1-q) J_0$,
and denoting by $|cm\rangle$ a simultaneous eigenvector of the commuting Hermitian
operators $C$ and $J_0$, associated with the eigenvalues $c$ and $m$ respectively,
it can be proved~\cite{delbecq1} that $J_+^n|cm\rangle$ (resp.~$J_-^n|cm\rangle$),
$n\in\N^+$, is either the null vector or a simultaneous eigenvector of $C$ and $J_0$,
corresponding to the eigenvalues $c$ and $mq^{-n}- \left(1-q^{-n}\right) / (1-q)$
(resp.~$mq^n- \left(1-q^n\right) / (1-q)$). Hence, the spectrum of $J_0$ is
exponential, instead of linear as for the PRA's.\par
%
%
Moreover, if the starting $m$ value belongs to the interval $((q-1)^{-1},+\infty)$
(resp.~$(-\infty,(q-1)^{-1})$), then all the $J_0$ eigenvalues obtained by successive
applications of $J_+$ or $J_-$ upon $|cm\rangle\ne0$ belong to the same interval
and $J_+$ (resp.\ $J_-$) is a raising generator, whereas if $m=(q-1)^{-1}$, then
neither $J_+$ nor $J_-$ change the $J_0$ eigenvalue. The unirreps therefore
separate into two classes according to whether the eigenvalues of $J_0$ are
contained in the interval $\bigl(-\infty, (q-1)^{-1}\bigr)$, or in the interval
$\bigl((q-1)^{-1}, +\infty\bigr)$.\par
%
%
In general, they may fall into one out of four categories: (i) infinite-dimensional
unirreps with a lower bound $-j$, (ii) infinite-dimensional unirreps with an upper
bound $J$, (iii) infinite-dimensional unirreps with neither lower nor upper bounds,
and (iv) finite-dimensional unirreps with both lower and upper bounds, $-j$ and $J$
(where in general $j\ne J$). In addition, there is a trivial one-dimensional unirrep
corresponding to $m=(q-1)^{-1}$.\par
%
%
The exponential character of the $J_0$~spectrum in DQA~representations may be of
interest in various physical problems, wherein such spectra have been encountered,
such as alternative Hamiltonian quantizations, exactly solvable potentials,
$q$-deformed supersymmetric quantum mechanics, and $q$-deformed interacting
boson models~\cite{fairlie}.\par
%
%
\section{The algebra \alg}
Another example of DQA, for which the function~$G(J_0)$ is linear, has been
recently constructed~\cite{bonatsos}. Contrary to those considered in the previous
section, this algebra, denoted by \alg, is defined in terms of functions~$F(J_0)$
and~$H(J_0)$ that are not polynomials, but infinite series in~$J_0$,
\begin{eqnarray}
  F(J_0) & = & - \frac{(G(J_0))^2 - (G(J_0))^{-2}}{q - q^{-1}}, \nonumber \\
  H(J_0) & = & \frac{q^{-1} (G(J_0))^2 + q (G(J_0))^{-2} - q - q^{-1}}{(q - q^{-1})^2},
  \label{eq:alg}
\end{eqnarray}
with $G(J_0) = 1 + (1-q) J_0$. Since the transformation $q \to q^{-1}$, $J_0 \to -q
J_0$, $J_{\pm} \to J_{\mp}$ is an automorphism of \alg, the parameter values
may be restricted to the range $0<q<1$.\par
%
%
\alg\ can be obtained from \su, a special case of PRA, by using a two-valued map
$P_{\delta} : \su \to \alg$, $\delta = \pm1$, defined by
\begin{equation}
  J_0 = p_{\delta}(j_0), \qquad J_+ = j_+, \qquad J_- = j_-, \qquad \mbox{where}
  \qquad p_{\delta}(z) \equiv \frac{1 - \delta q^{-z}}{q-1}.     \label{eq:map}
\end{equation}
Such a generator map is well defined: it can indeed be easily checked that if $j_0$,
$j_+$, $j_-$ satisfy the \su\ commutation relations, then $J_0$, $J_+$, $J_-$, given
in~(\ref{eq:map}), fulfil those of~\alg. The functions $p_{\delta}(z)$, $\delta =
\pm1$, defined in~(\ref{eq:map}), are entire and invertible functions, with $g(z) =
p_{\delta}^{-1}(z) = \ln\bigl(G^2(z)\bigr) \big/ \ln\bigl(q^{-2}\bigr)$. If $z \in \R$,
the range  of~$p_{\delta}$ (and consequently the domain of $p_{\delta}^{-1}$) is the
interval $\bigl(-\infty, (q-1)^{-1}\bigr)$ or $\bigl((q-1)^{-1}, +\infty\bigr)$
according to whether $\delta = -1$ or $\delta = +1$. The function~$g(z)$ is
well-behaved everywhere on~\R, except in the neighbourhood of the point $z =
(q-1)^{-1}$.\par
%
%
It should be stressed that the use of~$P_{\delta}$, $\delta = \pm1$, implies an
extension of the well-known deforming functional technique~\cite{curtright} for
two reasons: first because here a map between two deformed algebras, \su\ and \alg,
is considered instead of a map between a Lie algebra and a deformed one, as in the
original method; and second because use is made of a two-valued functional, whose
inverse is singular, instead of a single-valued one.\par
%
%
It can be easily shown~\cite{bonatsos} that \alg\ has no infinite-dimensional
unirrep, but has, for any $N = 0$, 1, 2,~$\ldots$, two ($N+1$)-dimensional unirreps,
which may be distinguished by $\delta = \pm1$. The corresponding spectrum of $J_0$
is given by $m^{\delta}  = \left(1 - \delta q^{-(N-2n)/2}\right) \big/ (q-1)$, $n = 0$,
1, $\ldots$,~$N$, with maximum and minimum eigenvalues $J^{\delta} =   \left(1
-\delta q^{-\delta N/2}\right) \big/ (q-1)$, and $-j^{\delta} =
\left(1 -\delta q^{\delta N/2}\right) \big/ (q-1)$ respectively. The unirrep
specified by $J^+$ (resp.~$J^-$) is entirely contained in the interval
$\bigl((q-1)^{-1}, +\infty \bigr)$ (resp.~$\bigl(-\infty, (q-1)^{-1}\bigr)$). For both
unirreps, the eigenvalue of the Casimir operator is given by $\langle C\rangle
=H(\gamma^\delta)$, where $\gamma^{\delta}= \left(1- \delta q^{-N/2}\right)
\big/ (q-1)$.\par
%
%
In the carrier space $V^{J^{\delta}}$ of the unirrep characterized by~$J^{\delta}$,
whose basis vectors are specified by the values of $J^{\delta}$ and~$m^{\delta}$,
the \alg\ generators are represented by some linear operators
$\Phi^{J^{\delta}}(A)$, $A \in \alg$, defined by
\begin{eqnarray}
   \Phi^{J^{\delta}}\left(J_0\right) \left\vert J^\delta, m^\delta \right> &=&
       m^\delta \left\vert J^\delta, m^\delta \right>
       = \left(\case{N}{2} - n\right) \left\vert J^\delta, m^\delta \right>,
       \nonumber\\[0.1cm]
   \Phi^{J^{\delta}}\left(J_-\right) \left\vert J^\delta, m^\delta \right> &=&
       \sqrt{ H(\gamma^{\delta}) - H(q m^\delta-1)}
       \left\vert J^\delta, qm^\delta-1 \right> \nonumber\\[0.1cm]
       & = & \sqrt{[n+1]_q [N-n]_q} \, \left\vert J^\delta, qm^\delta-1 \right>,
       \nonumber\\[0.1cm]
   \Phi^{J^{\delta}}\left(J_+\right) \left\vert J^\delta, m^\delta \right> &=&
       \sqrt{ H(\gamma^{\delta}) - H( m^\delta)}
       \left\vert J^\delta, q^{-1}(m^\delta +1) \right> \nonumber\\[0.1cm]
       & = & \sqrt{[n]_q [N-n+1]_q} \, \left\vert J^\delta, q^{-1}(m^\delta +1) \right>.  
       \label{eq:A-modules}
\end{eqnarray}
\par
%
%
The generator mapping~$P_{\delta}$ can be used to transfer the quasitriangular
Hopf structure of~\su\ to \alg~\cite{bonatsos}. One gets in this way a double
quasitriangular Hopf structure, with comultiplication, counit, antipode maps, and
universal \cR-matrix given by
\begin{eqnarray}
    \Delta_\delta \left(J_0 \right) &=& (q-1)^{-1} \left(1\otimes 1 - \delta 
        G(J_0)\otimes G(J_0) \right), \nonumber \\[0.1cm]
    \Delta_\delta \left( J_\pm \right) &=& \delta\left( J_{\pm} \otimes \left(
        G(J_0) \right)^{-1}  + G(J_0)\otimes J_{\pm} \right), \nonumber\\[0.1cm]
    \epsilon_\delta(J_0) &=& (1 - \delta) (q - 1)^{-1}, \quad 
        \epsilon_\delta (J_\pm) =0, \nonumber\\[0.1cm] 
    S_{\delta}(J_0) &=& - J_0 ( G(J_0) )^{-1}, \quad S_{\delta}(J_+) = - q J_+,
        \quad S_{\delta}(J_-) = - q^{-1} J_-, \nonumber \\[0.1cm]
     \cR^{\delta} &=& q^{2\log_q(\delta G(J_0)) \otimes  \log_q(\delta G(J_0))}
        \sum_{n=0}^{\infty} \frac{(1-q^{-2})^n}{[n]_q!}\, q^{n(n-1)/2} \nonumber
        \\[0.1cm]
     & & \mbox{} \times \left((G(J_0))^{-1} J_+ \otimes G(J_0) J_-\right)^n,  
        \label{eq:Com-def}
\end{eqnarray}
respectively. Both $\left(\Delta_+, \epsilon_+, S_+, \cR^+\right)$, and
$\left(\Delta_-,\epsilon_-, S_-,\cR^-\right)$ satisfy the Hopf and
quasitriangularity axioms, but the former are only valid for the representations of
$\alg$  with eigenvalues of $J_0$ in the interval $\bigl((q-1)^{-1}, +\infty \bigr)$,
whereas the latter act in $\bigl(-\infty, (q-1)^{-1}\bigr)$.\par
%
%
\section{Two-colour quasitriangular Hopf structure of \alg}

The double Hopf structure considered in the previous section allows one to couple
any two \alg\ unirreps  characterized by~$J_1^+$ and~$J_2^+$ (resp.~$J_1^-$
and~$J_2^-$), and with respective carrier spaces~$V^{J_1^+}$ and~$V^{J_2^+}$
(resp.~$V^{J_1^-}$ and~$V^{J_2^-}$), to obtain a reducible representation of the
same in $V^{J_1^+} \otimes V^{J_2^+}$ (resp.\ $V^{J_1^-} \otimes V^{J_2^-}$). No
coupling of two unirreps of the types~$J_1^+$ and~$J_2^-$, or~$J_1^-$ and~$J_2^+$,
is however possible.\par
%
%
To allow such types of couplings, it is necessary to extend the double Hopf
structure of \alg~\cite{bonatsos}. This can be accomplished by considering the
`transmutation' operators $T^{J^{\delta}}: V^{J^{\delta}} \to V^{J^{-\delta}}$, which
change the basis states of an ($N+1$)-dimensional unirrep, characterized
by~$J^{\delta}$, into those of its partner with the same dimension, specified
by~$J^{-\delta}$, i.e., $T^{J^{\delta}} \left\vert J^\delta, m^\delta \right> =
\left\vert J^{-\delta}, m^{-\delta}\right>$. By applying $T^{J^{\delta}}$ on both sides
of~(\ref{eq:A-modules}), we obtain that for any \alg\ generator~$A$,
$T^{J^{\delta}} \Phi^{J^{\delta}}(A)\, T^{J^{-\delta}} = \Phi^{J^{-\delta}}(\sigma(A))$,
where $\sigma: \alg \to \alg$, defined by $\sigma(J_0) = 2 (q-1)^{-1} - J_0$, and
$\sigma(J_{\pm}) = J_{\pm}$, is an involutive automorphism of the algebra \alg. This
clearly shows that at the algebra level, the operator~$\sigma$ is responsible for
the transmutation. Defining now $\sigma_{\delta}: \alg \to \alg$, $\delta = \pm1$, as
$\sigma_+ = \id$, and $\sigma_- = \sigma$, we note that the basic
mapping~$P_{\delta}$, defined in~(\ref{eq:map}), satisfies the relation 
$\sigma_{\zeta\eta} \circ P_{\eta} = P_{\zeta}$, where $\zeta$, $\eta=\pm 1$.\par
%
%
The comultiplication and antipode maps, as well as the double \cR-matrix of
equation~(\ref{eq:Com-def}) can be extended by setting
\begin{eqnarray}
  \Delta^{\zeta,\eta}_{\delta}(A)  & = & \left(\sigma_{\zeta \delta} \otimes
         \sigma_{\eta \delta}\right) \circ \Delta_{\delta}(A), \qquad
         S^{\zeta}_{\delta}(A) = \sigma_{\zeta \delta} \circ S_{\delta}(A), \nonumber
         \\[0.1cm]
  \cR^{\zeta,\eta} & = & (\sigma_{\zeta\delta} \otimes \sigma_{\eta\delta})
         \left(\cR^{\delta}\right),      \label{eq:extcom-def}
\end{eqnarray}
where $\zeta$, $\eta$, $\delta = \pm1$, while the counit map $\epsilon_{\delta}$,
defined in the same equation, is left unchanged. The results can be written as
\begin{eqnarray}
    \Delta^{\zeta,\eta}_\delta \left(J_0 \right) &=& (q-1)^{-1} \left(1\otimes 1
        - \delta \zeta \eta G(J_0)\otimes G(J_0) \right), \nonumber \\[0.1cm]
    \Delta^{\zeta,\eta}_\delta \left( J_\pm \right) &=& \eta J_{\pm} \otimes
        \left(G(J_0) \right)^{-1}  + \zeta G(J_0)\otimes J_{\pm},
        \nonumber\\[0.1cm]
    S^{\zeta}_{\delta}(J_0) &=& (q-1)^{-1} \left(1 - \zeta \delta
        \left(G(J_0)\right)^{-1}\right), \nonumber \\[0.1cm] 
    S^{\zeta}_{\delta}(J_{\pm}) & = & - q^{\pm1} J_{\pm}, \nonumber \\[0.1cm]
    \cR^{\zeta,\eta} & = & q^{2\log_q(\zeta G(J_0)) \otimes  \log_q(\eta G(J_0))}
        \nonumber \\[0.1cm] 
    & & \mbox{} \times \sum_{n=0}^{\infty} \frac{(1-q^{-2})^n}{[n]_q!}\, q^{n(n-1)/2}
       \left((\zeta G(J_0))^{-1} J_+ \otimes \eta G(J_0) J_-\right)^n.      
       \label{eq:extcom}
\end{eqnarray}
\par
%
%
It can be easily shown~\cite{bonatsos} that the generalized comultiplication, counit,
and antipode maps, $\Delta^{\zeta,\eta}_{\delta}$, $\epsilon_{\delta}$,
$S^{\zeta}_{\delta}$, defined in~(\ref{eq:Com-def}) and~(\ref{eq:extcom}),
transform under $\sigma_{\delta}$ as 
\begin{equation}
  \left(\sigma_{\mu \zeta} \otimes \sigma_{\nu \eta}\right) \circ
      \Delta^{\zeta,\eta}_{\delta} = \Delta^{\mu,\nu}_{\rho} \circ
      \sigma_{\rho\delta}, \qquad
  \epsilon_{\delta} \circ \sigma_{\delta \zeta} = \epsilon_{\zeta}, \qquad
  \sigma_{\zeta \eta} \circ S^{\eta}_{\delta} = S^{\zeta}_{\mu} \circ
      \sigma_{\mu \delta},         \label{eq:extcom-transf3}
\end{equation}
and satisfy the following generalized coassociativity, counit, and
antipode axioms, 
\begin{eqnarray}
  \left(\Delta^{\zeta,\eta}_{\mu} \otimes \id\right) \circ \Delta^{\mu,\nu}_{\delta}
        (A) & = & \left(\id \otimes \Delta^{\eta,\nu}_{\rho}\right) \circ
        \Delta^{\zeta,\rho}_{\delta}(A), \nonumber \\[0.1cm]
  \left(\epsilon_{\zeta} \otimes \sigma_{\eta \delta}\right) \circ
        \Delta^{\zeta,\eta}_{\delta}(A) & = & \left(\sigma_{\zeta \delta} \otimes
        \epsilon_{\eta}\right) \circ \Delta^{\zeta,\eta}_{\delta}(A) = A, 
        \nonumber\\[0.1cm]
  m \circ \left(S^{\mu}_{\zeta} \otimes \sigma_{\mu \eta}\right) \circ
        \Delta^{\zeta,\eta}_{\delta}(A) & = & m \circ \left(\sigma_{\mu \zeta}
        \otimes S^{\mu}_{\eta}\right) \circ \Delta^{\zeta,\eta}_{\delta}(A)
        \nonumber \\[0.1cm] 
        & = & \iota \circ \epsilon_{\delta}(A),      \label{eq:genHopf}
\end{eqnarray}
where $A$ denotes any element of \alg, $m$ and $\iota$ are the multiplication and
unit maps of \alg, $\delta$, $\zeta$, $\eta$, $\mu$, $\nu$, $\rho$ take any values in
the set $\{-1,+1\}$, and no summation over repeated indices is implied. Moreover,
$\Delta^{\zeta,\eta}_{\delta}$ and $\epsilon_{\delta}$ are algebra homomorphisms,
while $S^{\zeta}_{\delta}$ is both an algebra and a coalgebra antihomomorphism.\par
%
%
By using the generalized coproduct $\Delta^{\zeta,\eta}_{\delta}$, it is now
possible to couple any ($N_1+1$)- and ($N_2+1$)-dimensional unirreps of \alg,
specified by $J_1^{\zeta}$ and $J_2^{\eta}$ respectively, to provide two reducible
representations in $V^{J_1^{\zeta}} \otimes V^{J_2^{\eta}}$, which are characterized
by $\delta = +1$ and $\delta = -1$, respectively. They can be decomposed into a
direct sum of ($N+1$)-dimensional unirreps, specified by~$J^{\delta}$, by using some
Wigner coefficients $\bigl\langle J^{\zeta}_1\, m^{\zeta}_1,  J^{\eta}_2\,
m^{\eta}_2 \big| J^{\delta}\, m^{\delta}\bigr\rangle_{DQ}$, given in terms
of \su\ Wigner coefficients by the relation
\begin{equation}
   \bigl\langle J^{\zeta}_1\, m^{\zeta}_1,  J^{\eta}_2\, m^{\eta}_2 
   \big| J^{\delta}\, m^{\delta} \bigr\rangle_{DQ} = 
   \bigl\langle \case{N_1}{2}\; \case{N_1}{2} - n_1, 
   \case{N_2}{2}\; \case{N_2}{2} - n_2 \big| \case{N}{2}\; \case{N}{2} - n
   \bigr\rangle_q.     \label{eq:Wigner}
\end{equation}
The carrier space of the unirrep~$J^{\delta}$ in $V^{J_1^{\zeta}} \otimes
V^{J_2^{\eta}}$ is therefore spanned by the states
\begin{equation}
  \bigl| J_1^{\zeta} J_2^{\eta} J^{\delta} m^{\delta} \bigr\rangle = 
  \sum_{m_1^{\zeta}, m_2^{\eta}} \bigl\langle J^{\zeta}_1\, m^{\zeta}_1,  
  J^{\eta}_2\, m^{\eta}_2 \big| J^{\delta}\, m^{\delta}
  \bigr\rangle_{DQ} \bigl|J_1^{\zeta}, m_1^{\zeta}\bigr\rangle \otimes
  \bigl|J_2^{\eta}, m_2^{\eta} \bigr\rangle.   \label{eq:coupled}  
\end{equation}
\par
%
%
Turning now to the generalized universal \cR-matrix defined
in~(\ref{eq:extcom-def}) or~(\ref{eq:extcom}), it can be easily
shown~\cite{bonatsos} that its four pieces $\cR^{\zeta,\eta}$, $\zeta$, $\eta =
\pm1$, are invertible and satisfy the properties
\begin{eqnarray}
  (\sigma_{\mu\zeta} \otimes \sigma_{\nu\eta}) \left(\cR^{\zeta,\eta}\right) & =
          & \cR^{\mu,\nu}, \quad \tau \circ \Delta^{\eta,\zeta}_{\delta}(A) = 
          \cR^{\zeta,\eta} \Delta^{\zeta,\eta}_{\delta}(A)
          \left(\cR^{\zeta,\eta}\right)^{-1}, \nonumber \\
  \left(\Delta^{\lambda,\mu}_{\zeta} \otimes \sigma_{\nu\eta}\right) 
          \left(\cR^{\zeta,\eta}\right) & = & \cR^{\lambda,\nu}_{13}\,
          \cR^{\mu,\nu}_{23}, \quad \left(\sigma_{\lambda\zeta} \otimes
          \Delta^{\mu,\nu}_{\eta}\right) \left(\cR^{\zeta,\eta}\right) = 
          \cR^{\lambda,\nu}_{13}\, \cR^{\lambda,\mu}_{12},    \label{eq:genR-prop1}  
\end{eqnarray} 
for any $A \in \alg$. From these results, or more simply from the corresponding
properties fulfilled by~$\cR^{\delta}$, one obtains the relations 
\begin{eqnarray}
  \cR^{\zeta,\eta}_{12}\, \cR^{\zeta,\mu}_{13}\, \cR^{\eta,\mu}_{23} & = &
         \cR^{\eta,\mu}_{23}\, \cR^{\zeta,\mu}_{13}\, \cR^{\zeta,\eta}_{12},     
         \nonumber \\
  (\epsilon_{\zeta} \otimes \id) \left(\cR^{\zeta,\eta}\right) & = & (\id \otimes
         \epsilon_{\eta}) \left(\cR^{\zeta,\eta}\right) = 1, \nonumber \\
  (S^{\lambda}_{\zeta} \otimes \sigma_{\mu\eta}) \left(\cR^{\zeta,\eta}\right) & = &
         \left(\sigma_{\lambda\zeta} \otimes (S^{\mu}_{\eta})^{-1}\right)
         \left(\cR^{\zeta,\eta}\right) = \left(\cR^{\lambda,\mu}\right)^{-1}.
         \label{eq:genR-prop2}
\end{eqnarray}
\par
%
%
The first relation in~(\ref{eq:genR-prop2}) shows that the generalized universal
\cR-matrix is a solution of the coloured YBE~\cite{murakami}, where the `colour'
parameters $\zeta$, $\eta$,~$\mu$ take discrete values in the set $\{-1,+1\}$. We
may therefore call
$\Bigl(\alg,+,m,\iota,\Delta^{\zeta,\eta}_{\delta},\epsilon_{\delta},
S^{\zeta}_{\delta},\cR^{\zeta,\eta};\C\Bigr)$ a two-colour quasitriangular Hopf
algebra over \C~\cite{bonatsos}. As will be shown elsewhere~\cite{cq}, this type of
algebraic structure admits generalizations, which will be referred to as coloured
quasitriangular Hopf algebras.
\par
%
%

\end{document}